\documentclass[12pt,oneside,reqno]{amsart}
\usepackage{mathrsfs}
\usepackage{graphics}
\usepackage{amssymb}
\usepackage{enumerate}
\newcommand{\dif}{\mathrm{d}}

\newcommand{\be}{\begin{eqnarray}}
\newcommand{\ee}{\end{eqnarray}}
\newcommand{\ce}{\begin{eqnarray*}}
\newcommand{\de}{\end{eqnarray*}}
\newtheorem{theorem}{Theorem}[section]
\newtheorem{lemma}[theorem]{Lemma}
\newtheorem{remark}[theorem]{Remark}
\newtheorem{definition}[theorem]{Definition}
\newtheorem{proposition}[theorem]{Proposition}
\newtheorem{Example}[theorem]{Example}
\newtheorem{corollary}[theorem]{Corollary}
\def\e{\varepsilon}
\def\s{\sigma}

\def\[{{\Big[}}
\def\]{{\Big]}}
\def\<{{\langle}}
\def\>{{\rangle}}
\def\({{\Big(}}
\def\){{\Big)}}

\def\no{\nonumber}
\def\bt{\begin{theorem}}
\def\et{\end{theorem}}
\def\bl{\begin{lemma}}
\def\el{\end{lemma}}
\def\br{\begin{remark}}
\def\er{\end{remark}}
\def\bx{\begin{Example}}
\def\ex{\end{Example}}
\def\bd{\begin{definition}}
\def\ed{\end{definition}}
\def\bp{\begin{proposition}}
\def\ep{\end{proposition}}
\def\bc{\begin{corollary}}
\def\ec{\end{corollary}}

\def\cB{{\mathcal B}}

\def\cD{{\mathcal D}}

\def\cK{{\mathcal K}}
\def\cL{{\mathcal L}}
\def\cM{{\mathcal M}}

\def\cP{{\mathcal P}}

\def\mE{{\mathbb E}}

\def\mL{{\mathbb L}}

\def\mN{{\mathbb N}}

\def\mP{{\mathbb P}}

\def\mR{{\mathbb R}}

\def\mW{{\mathbb W}}

\def\sB{{\mathscr B}}
\def\sC{{\mathscr C}}

\def\sF{{\mathscr F}}
\def\sG{{\mathscr G}}

\def\sL{{\mathscr L}}

\def\sS{{\mathscr S}}

\def\geq{\geqslant}
\def\leq{\leqslant}

\begin{document}

\allowdisplaybreaks

\title{Uniqueness and superposition of the space-distribution dependent Zakai equations*}

\author{Meiqi Liu$^1$ and Huijie Qiao$^{1,2}$}

\dedicatory{1. Department of Mathematics,
Southeast University,\\
Nanjing, Jiangsu 211189, P.R.China\\
2. Department of Mathematics, University of Illinois at
Urbana-Champaign\\
Urbana, IL 61801, USA}



\thanks{*This work was partly supported by NSF of China (No. 11001051, 11371352, 12071071) and China Scholarship Council under Grant No. 201906095034.}

\thanks{Corresponding author: Huijie Qiao, hjqiaogean@seu.edu.cn}


\date{}

\begin{abstract}
The work concerns the space-distribution dependent Zakai equations from nonlinear filtering problems of McKean-Vlasov stochastic differential equations with correlated noises. First of all, we establish the space-distribution dependent Kushner-Stratonovich equations and the space-distribution dependent Zakai equations. Then, the pathwise uniqueness of their strong solutions is shown. Finally, we prove a superposition principle between the space-distribution dependent Zakai equations and space-distribution dependent Fokker-Planck equations. As a by-product, we give some conditions under which space-distribution dependent Fokker-Planck equations have weak solutions.
\end{abstract}

\maketitle \rm

\section{Introduction}

McKean-Vlasov (distribution-dependent or mean-field) stochastic differential equations (SDEs for short) describe the evolution rules of particle systems perturbed by noises. The difference between McKean-Vlasov SDEs and general SDEs is that the former  depend on the positions and probability distributions of these particles. Therefore, McKean-Vlasov SDEs are widely applied in many fields, such as biology, game theory and control theory. Moreover, more and more results about McKean-Vlasov SDEs appear. We mention some results associated with our work. Ding and Qiao \cite{dq1, dq2} investigated the well-posedness and stability of weak solutions for McKean-Vlasov SDEs under non-Lipschitz conditions. Lacker, Shkolnikov and Zhang \cite{lsz} studied superposition principles for conditional McKean-Vlasov equations. Ren and Wang \cite{rw} proved that additive functionals of McKean-Vlasov SDEs have path-independence.

Nonlinear filtering problems are to extract some useful information of unobservable phenomenon from observable ones, and estimate and predict them (c.f. \cite{Crisan, ko, kurtz-xiong, q1, q11, q2, qd, roge, blr}). Thus, nonlinear filtering theory plays an important role in many areas including stochastic control, financial modeling, speech and image processing, and Bayesian networks. Although McKean-Vlasov SDEs have widespread applications, the result about nonlinear filtering problems of McKean-Vlasov SDEs is seldom. Only Sen and Caines \cite{sc, sc1} studied nonlinear filtering problems of McKean-Vlasov SDEs with independent noises. 

In the paper, we focus on nonlinear filtering problems of McKean-Vlasov SDEs with correlated noises. We explain them in detail. Fix $T>0$. Let $(\Omega,\sF,\{\sF_t\}_{t\in[0,T]},\mP)$ be a complete filtered probability space and $\{W_{t},t\geq0\}$, $\{V_{t},t\geq0\}$ be $d$-dimensional and $m$-dimensional standard Brownian motions defined on $(\Omega,\sF,\{\sF_t\}_{t\in[0,T]},\mP)$, respectively. Moreover, $W_{\cdot}$ and $V_{\cdot}$ are mutually independent. Consider the following McKean-Vlasov signal-observation system $(X_t, Y_t)$ on $\mR^n\times\mR^m$:
\be\left\{\begin{array}{l}
\dif X_t=b_1(t,X_t,\sL^{\mP}_{X_t})\dif t+\sigma_0(t,X_t,\sL^{\mP}_{X_t})\dif W_t+\sigma_1(t,X_t,\sL^{\mP}_{X_t})\dif V_t,\\
\dif Y_t=b_2(t,X_t,\sL^{\mP}_{X_t},Y_t)\dif t+\sigma_2(t,Y_t)\dif V_t, \quad 0\leq t\leq T,
 \end{array}
\right. 
\label{Eq1} 
\ee
where $\sL^{\mP}_{X_t}$ denotes the distribution of $X_t$ under the probability measure $\mP$, and these coefficients $b_1: [0,T] \times \mR^n \times \cP_2(\mR^n) \mapsto \mR^n$, $\s_0:  [0,T] \times \mR^n \times \cP_2(\mR^n)\mapsto\mR^{n\times d}$, $\s_1:  [0,T] \times \mR^n \times \cP_2(\mR^n)\mapsto\mR^{n\times m}$, $b_2: [0,T] \times \mR^n \times \cP_2(\mR^n)\times\mR^m \mapsto \mR^m$ and $\s_2:  [0,T] \times \mR^m\mapsto\mR^{m\times m}$ are Borel measurable. The initial value $X_0$ is assumed to be a $p$-order ($p>2$) integrable random variable independent of $Y_0, W_{\cdot}, V_{\cdot}$. The system (\ref{Eq1}) is called a model with a correlated noise. Then we deduce the space-distribution dependent Kushner-Stratonovich equation and the space-distribution dependent Zakai equation about the system (\ref{Eq1}). Next, we view  the space-distribution dependent Kushner-Stratonovich equation and the space-distribution dependent Zakai equation as two SDEs and define their strong solutions. Moreover, we prove that strong solutions of the space-distribution dependent Kushner-Stratonovich equation and the space-distribution dependent Zakai equation both have the pathwise uniqueness. Finally, we define weak solutions of the space-distribution dependent Zakai equation and a space-distribution dependent Fokker-Planck equation, and set up a correspondence between weak solutions of the space-distribution dependent Zakai equation and weak solutions of a space-distribution dependent Fokker-Planck equation.

Moreover, our methods can be applied to study nonlinear filtering problems of McKean-Vlasov SDEs with correlated sensor noises. Concretely speaking, consider the following signal-observation system $(\check{X}_t, \check{Y}_t)$ on $\mR^n\times\mR^m$:
\be\left\{\begin{array}{l}
\dif \check{X}_t=\check{b}_1(t,\check{X}_t, \sL^{\mP}_{\check{X}_t})\dif t+\check{\sigma}_1(t,\check{X}_t, \sL^{\mP}_{\check{X}_t})\dif V_t,\\
\dif \check{Y}_t=\check{b}_2(t,\check{X}_t, \sL^{\mP}_{\check{X}_t},\check{Y}_t)\dif t+\check{\sigma}_2\dif W_t+\check{\sigma}_3\dif V_t, 
 \end{array}
\right. 
0\leq t\leq T,
\label{Eq0} 
\ee
where the initial value $\check{X}_0$ is assumed to be a $p$-order ($p>2$) integrable random variable independent of $\check{Y}_0, W_{\cdot}, V_{\cdot}$. The mappings $\check{b}_1:[0,T]\times\mR^n\times \cP_2(\mR^n) \mapsto\mR^n$, $\check{\sigma}_1:[0,T]\times\mR^n\times \cP_2(\mR^n) \mapsto\mR^{n\times m}$ and $\check{b}_2:[0,T]\times\mR^n\times \cP_2(\mR^n)\times\mR^m\mapsto\mR^m$ are all Borel measurable. $\check{\sigma}_2, \check{\sigma}_3$ are $m\times d$ and $m\times m$ real matrices, respectively. By the way similar to that for the system (\ref{Eq1}), we can also establish another space-distribution dependent Zakai equation, and study its pathwise uniqueness and superposition principles.

The paper is arranged as follows. In Section \ref{pre}, we introduce notation and $L$-derivative for functions on $\cP_2(\mR^n)$ used in the sequel. After this, we introduce nonlinear filtering problems for McKean-Vlasov signal-observation systems with correlated noises,  and derive the space-distribution dependent Kushner-Stratonovich equations and the space-distribution dependent Zakai equations. In Section \ref{unzaks}, the pathwise uniqueness for strong solutions to the space-distribution dependent Kushner-Stratonovich equations and the space-distribution dependent Zakai equations is shown. We place a superposition principle for the space-distribution dependent Zakai equation in Section \ref{super}. Finally, in Section \ref{con} we summarize our results and apply our methods to the system (\ref{Eq0}).

The following convention will be used throughout the paper: $C$, with or without indices, will denote
different positive constants whose values may change from one place to another.

\section{Preliminary}\label{pre}

In the section, we introduce notation and $L$-derivative for functions on $\cP_2(\mR^n)$.

\subsection{Notation}\label{nn}

In the subsection, we introduce notation used in the sequel. 

For convenience, we shall use $\mid\cdot\mid$ and $\parallel\cdot\parallel$  for norms of vectors and matrices, respectively. Let $A^*$ denote the transpose of the matrix $A$.

Let $\sB(\mR^n)$ be the Borel $\sigma$-field on $\mR^n$. Let $\cB_b(\mR^n)$ denote the set of all real-valued uniformly bounded $\mathscr{B}(\mR^n)$-measurable functions on $\mR^n$. $C^2(\mR^n)$ stands for the space of continuous functions on $\mR^n$ which have continuous partial derivatives of order up to $2$, and $C_b^2(\mR^n)$ stands for the subspace of $C^2(\mR^n)$, consisting of functions whose derivatives up to order 2 are bounded. $C_c^2(\mR^n)$ is the collection of all functions in $C^2(\mR^n)$ with compact supports and $C_c^\infty(\mR^n)$ denotes the collection of all real-valued $C^\infty$ functions of compact
supports.

Let $\cM(\mR^n)$ be the set of all bounded Borel measures defined on $\sB(\mR^n)$ carrying the usual topology of weak convergence. Let $\cP({\mR^n})$ be the space of all probability measures defined on $\sB(\mR^n)$ and $\cP_2(\mR^n)$ be the collection of all the probability measures $\mu$ on $\sB(\mR^n)$ satisfying
\ce
\|\mu\|_2^2:=\int_{\mR^n}\mid{x}\mid^2\,\mu(\dif x)<\infty.
\de
We put on $\cP_2(\mR^n)$ a topology induced by the
following $2$-Wasserstein metric:
\ce
\mW_2^2(\mu_1,\mu_2):=\inf_{\pi\in\sC(\mu_1, \mu_2)}\int_{\mR^n\times\mR^n}|x-y|^2\pi(\dif x, \dif y), \quad \mu_1, \mu_2\in\cP_2(\mR^n),
 \de
where $\sC(\mu_1, \mu_2)$ denotes the set of  all the probability measures whose marginal distributions are $\mu_1, \mu_2$, respectively. It is known that $(\cP_2(\mR^n), \mW_2)$ is a Polish space (c.f. \cite{vc}).

\subsection{$L$-derivative for functions on $\cP_2(\mR^n)$}\label{lde} In the subsection we recall the definition of $L$-derivative for functions on $\cP_2(\mR^n)$. The definition was first introduced by Lions (c.f. \cite{Lion}). Moreover, he used some abstract probability spaces to describe the $L$-derivatives. Here, for the convenience to understand the definition, we apply a straight way to state it (c.f. \cite{rw}). Let $I$ be the identity map on $\mR^n$. For $\mu\in\cP_2(\mR^n)$ and $\phi\in L^2(\mR^n, \sB(\mR^n), \mu;\mR^n)$, $<\mu,\phi>:=\int_{\mR^n}\phi(x)\mu(\dif x)$. Moreover, by simple calculation, it holds that $\mu\circ(I+\phi)^{-1}\in\cP_2(\mR^n)$.

\bd\label{lderi}
(i) A function $f: \cP_2(\mR^n)\mapsto\mR$ is called L-differentiable at $\mu\in\cP_2(\mR^n)$, if the functional 
$$
L^2(\mR^n, \sB(\mR^n), \mu;\mR^n)\ni\phi\mapsto f(\mu\circ(I+\phi)^{-1})
$$
is Fr\'echet differentiable at $\phi=0$; that is, there exists a unique $\gamma\in L^2(\mR^n, \sB(\mR^n), \mu;\mR^n)$ such that 
$$
\lim\limits_{<\mu,|\phi|^2>\rightarrow 0}\frac{f(\mu\circ(I+\phi)^{-1})-f(\mu)-<\mu,\gamma\cdot\phi>}{\sqrt{<\mu,|\phi|^2>}}=0.
$$
In the case, we denote $\partial_{\mu}f(\mu)=\gamma$ and call it the L-derivative of $f$ at $\mu$.

(ii) A function $f: \cP_2(\mR^n)\mapsto\mR$ is called L-differentiable on $\cP_2(\mR^n)$ if  L-derivative $\partial_{\mu}f(\mu)$ exists for all $\mu\in\cP_2(\mR^n)$.

(iii) By the same way, $\partial^2_\mu f(\mu)(y,y')$ for $y, y'\in\mR^n$ can be defined.
\ed

Next, we introduce some related spaces.

\bd\label{space1}
 The function $f$ is said to be in $C^2(\cP_2(\mR^n))$, if $\partial_\mu f$ is continuous, for any $\mu\in \cP_2(\mR^n)$, $\partial_\mu f(\mu)(\cdot)$ is differentiable, and its derivative $\partial_y\partial_\mu f:\cP_2(\mR^n)\times\mR^n\mapsto\mR^n\otimes\mR^n$ is continuous, and for any $y\in\mR^n$, $\partial_\mu f(\cdot)(y)$ is differentiable, and its derivative $\partial^2_\mu f:\cP_2(\mR^n)\times\mR^n\times\mR^n\mapsto\mR^n\otimes\mR^n$ is continuous.
\ed

\bd\label{space2}
(i) The function $F: \mR^n\times\cP_2(\mR^n)\mapsto\mR$ is said to be in $C^{2,2}(\mR^n\times\cP_2(\mR^n))$, if $F(x,\mu)$ is $C^2$ in $x\in\mR^n$ and $\mu\in\cP_2(\mR^n)$ respectively, and its derivatives 
$$
\partial_x F(x,\mu), \partial^2_x F(x,\mu), \partial_\mu F(x,\mu)(y),  \partial_y\partial_\mu F(x,\mu)(y), \partial^2_\mu F(x,\mu)(y, y')
$$ 
are jointly continuous in the corresponding variable family $(x,\mu)$, $(x,\mu,y)$ or $(x,\mu,y, y')$. 

(ii) The function $F: \mR^n\times\cP_2(\mR^n)\mapsto\mR$ is said to be in $C_b^{2,2}(\mR^n\times\cP_2(\mR^n))$, if $F$ belongs to $C^{2,2}(\mR^n\times\cP_2(\mR^n))$ and is uniformly continuous with respect to $(x,\mu)$, and its derivatives and itself 
are bounded.

(iii) The function $F: \mR^n\times\cP_2(\mR^n)\mapsto\mR$ is said to be in $\sS(\mR^n\times\cP_2(\mR^n))$, if $F\in C^{2,2}(\mR^n\times\cP_2(\mR^n))$ and for any compact set $\cK\subset\mR^n\times\cP_2(\mR^n)$,
$$
\sup\limits_{(x,\mu)\in\cK}\int_{\mR^n}\left(\|\partial_y\partial_\mu F(x,\mu)(y)\|^2+|\partial_\mu F(x,\mu)(y)|^2\right)\mu(\dif y)<\infty.
$$

(iv) The function $\Phi: \mR^n\times\cP_2(\mR^n)\times\mR^m\mapsto\mR$ is said to be in $C^{2,2,2}(\mR^n\times\cP_2(\mR^n)\times\mR^m)$, if for $y\in\mR^m$, $\Phi(\cdot,\cdot,y)\in \sS(\mR^n\times\cP_2(\mR^n))$ and for $(x,\mu)\in\mR^n\times\cP_2(\mR^n)$, $\Phi(x,\mu,\cdot)\in C^2(\mR^m)$.
\ed

\section{Nonlinear filtering problems for McKean-Vlasov signal-observation systems with correlated noises}\label{nonfilter}

In this section, we introduce nonlinear filtering problems for McKean-Vlasov signal-observation systems with correlated noises,  and derive the space-distribution dependent Kushner-Stratonovich equations and the space-distribution dependent Zakai equations.

\subsection{The framework}\label{fram}

In the subsection, we introduce McKean-Vlasov signal-observation systems.

Consider the system (\ref{Eq1}), i.e. 
\ce\left\{\begin{array}{l}
\dif X_t=b_1(t,X_t,\sL^{\mP}_{X_t})\dif t+\sigma_0(t,X_t,\sL^{\mP}_{X_t})\dif W_t+\sigma_1(t,X_t,\sL^{\mP}_{X_t})\dif V_t,\\
\dif Y_t=b_2(t,X_t,\sL^{\mP}_{X_t},Y_t)\dif t+\sigma_2(t,Y_t)\dif V_t, \quad 0\leq t\leq T.
 \end{array}
\right. 
\de
We assume the following:
\begin{enumerate}[($\mathbf{H}^1_{b_1, \sigma_0, \sigma_1}$)] 
\item For $t\in[0,T]$ and $x_1, x_2\in\mR^n$, $\mu_1, \mu_2\in\cP_2(\mR^n)$,
\ce
&|b_1(t,x_1,\mu_1)-b_1(t,x_2,\mu_2)|\leq L_1(t)\(|x_1-x_2|\kappa_1(|x_1-x_2|)+\mW_2(\mu_1,\mu_2)\),\\
&\|\sigma_0(t,x_1,\mu_1)-\sigma_0(t,x_2,\mu_2)\|^2\leq L_1(t)\(|x_1-x_2|^{2}\kappa_2(|x_1-x_2|)+\mW_2^2(\mu_1,\mu_2)\),\\
&\|\sigma_1(t,x_1,\mu_1)-\sigma_1(t,x_2,\mu_2)\|^2\leq L_1(t)\(|x_1-x_2|^{2}\kappa_3(|x_1-x_2|)+\mW_2^2(\mu_1,\mu_2)\),
\de
where $L_1(t)>0$ is an increasing function and $\kappa_i$ is a positive continuous
function, bounded on $[1,\infty)$ and satisfies
\ce
\lim\limits_{x\downarrow0}\frac{\kappa_i(x)}{\log x^{-1}}<\infty, \quad i=1, 2, 3.
\de
\end{enumerate}
\begin{enumerate}[($\mathbf{H}^2_{b_1, \sigma_0,\sigma_1}$)]
\item For $t\in[0,T]$ and $x\in\mR^n$, $\mu\in\cP_2(\mR^n)$,
$$
|b_1(t,x,\mu)|^2+\|\sigma_0(t,x,\mu)\|^2+\|\sigma_1(t,x,\mu)\|^2\leq K_1(t)(1+|x|+\|\mu\|_2)^2,
$$
where $K_1(t)>0$ is an increasing function.
\end{enumerate}
\begin{enumerate}[($\mathbf{H}^1_{\sigma_2}$)] 
\item For $t\in[0,T]$ and $y_1, y_2\in\mR^m$,
\ce
\|\sigma_2(t,y_1)-\sigma_2(t,y_2)\|^2\leq L_2(t)|y_1-y_2|^2,
\de
where $L_2(t)>0$ is an increasing function.
\end{enumerate}
\begin{enumerate}[($\mathbf{H}^2_{b_2, \sigma_2}$)] 
\item For $t\in[0,T], y\in\mR^m$, $\sigma_2(t,y)$ is invertible, and
$$
|b_2(t,x,\mu, y)|\vee\|\sigma_2(t,0)\|\vee\|\sigma^{-1}_2(t,y)\|\leq K_2,~{for}~{all}~t\in[0,T], x\in\mR^n, \mu\in\cP_2(\mR^n), y\in\mR^m,
$$
where $K_2>0$ is a constant.
\end{enumerate}

Under the assumptions ($\mathbf{H}^1_{b_1, \sigma_0, \sigma_1}$), ($\mathbf{H}^2_{b_1, \sigma_0,\sigma_1}$), ($\mathbf{H}^1_{\sigma_2}$), ($\mathbf{H}^2_{b_2, \sigma_2}$),  by Theorem 3.1 in \cite{dq1}, it holds that the system (\ref{Eq1}) has a pathwise unique strong solution denoted as $(X_t,Y_t)$. Set
\ce
h(t,x,\mu,y):=\sigma_2^{-1}(t,y)b_2(t,x,\mu, y),
\de
\ce
\Gamma^{-1}_t:&=&\exp\bigg\{-\int_0^th^i(s,X_s,\sL^{\mP}_{X_s},Y_s)\dif V^i_s-\frac{1}{2}\int_0^t
\left|h(s,X_s,\sL^{\mP}_{X_s},Y_s)\right|^2\dif s\bigg\}.
\de
Here and hereafter, we use the convention that repeated indices imply summation. By ($\mathbf{H}^2_{b_2, \sigma_2}$), we know that 
\ce
\mE\left[\exp\left\{\int_0^T
\left|h(s,X_s,\sL^{\mP}_{X_s},Y_s)\right|^2\dif s\right\}\right]<\infty.
\de
Thus, the Novikov condition holds, and $\Gamma^{-1}_{\cdot}$ is an exponential martingale. Define a measure $\tilde{\mP}$ via
$$
\frac{\dif \tilde{\mP}}{\dif \mP}=\Gamma^{-1}_T,
$$
and under the measure $\tilde{\mP}$, 
\be\label{tilw}
\tilde{V}_t:=V_t+\int_0^t h(s,X_s,\sL^{\mP}_{X_s},Y_s)\dif s
\ee
is an $(\mathscr{F}_t)$-adapted Brownian motion. Moreover, the $\sigma$-algebra $\mathscr{F}_t^{Y}$ generated by $\{Y_s, 0\leq s\leq t\}$, can be characterized as
\ce
\mathscr{F}_t^{Y}=\mathscr{F}_t^{\tilde{V}}\vee\mathscr{F}_0^{Y},
\de
where $\mathscr{F}_t^{\tilde{V}}$ denotes the $\sigma$-algebra generated by $\{\tilde{V}_s, 0\leq s\leq t\}$. We augment $\mathscr{F}_t^{Y}$ in a usual sense and still denote the augmentation of $\mathscr{F}_t^{Y}$ as $\mathscr{F}_t^Y$. 

\subsection{The space-distribution dependent Kushner-Stratonovich equation}\label{kseq}

Set
\ce
<\Lambda_t,F>:=\mE[F(X_t,\sL^{\mP}_{X_t})|\mathscr{F}_t^Y], \quad F\in\cB_b(\mR^n\times\cP_2(\mR^n)),
\de
where $\cB_b(\mR^n\times\cP_2(\mR^n))$ denotes the set of all bounded measurable functions on $\mR^n\times\cP_2(\mR^n)$, and then $\Lambda_t$ is called the nonlinear filtering of $(X_t,\sL^{\mP}_{X_t})$ with respect to $\mathscr{F}_t^Y$. Moreover, the equation $\Lambda_{\cdot}$ satisfying is called the space-distribution dependent Kushner-Stratonovich equation. In order  to derive the space-distribution dependent Kushner-Stratonovich equation, we need the following result (c.f. Lemma 2.2 in \cite{q1}).

\bl\label{brmoposs}
Under the measure $\mP$, $\bar{V}_t:=\tilde{V}_t-\int_0^t<\Lambda_s,h(s,\cdot,\cdot,Y_s)>\dif s$ is an $(\mathscr{F}^Y_t)$-adapted Brownian motion.
\el

Now, it is the position to establish the space-distribution dependent Kushner-Stratonovich equation.

\bt  (The space-distribution dependent Kushner-Stratonovich equation)  \label{ks}\\
For $F\in\sS(\mR^n\times\cP_2(\mR^n))$, the space-distribution dependent Kushner-Stratonovich equation of the system (\ref{Eq1}) is given by
\be
<\Lambda_t,F>&=&<\Lambda_0,F>+\int_0^t<\Lambda_s,\mL_s F>\dif s+\int_0^t<\Lambda_s,\partial_{x_i} F\sigma^{ij}_1(s,\cdot,\cdot)>\dif  \bar{V}^j_s\no\\
&&+\int_0^t\left(<\Lambda_s,Fh^j(s,\cdot,\cdot,Y_s)>-<\Lambda_s,F><\Lambda_s,h^j(s,\cdot,\cdot,Y_s)>\right)\dif \bar{V}^j_s, \no\\
&&\qquad\qquad\qquad\qquad t\in[0,T], \label{kseq0}
\ee
where the operater $\mL_s$ is defined as
\be
(\mL_s F)(x,\mu)&=&\partial_{x_i}F(x,\mu)b^i_1(s,x,\mu)+\frac{1}{2}\partial_{x_ix_j}^2F(x,\mu)(\sigma_0\sigma_0^*)^{ij}(s,x,\mu)\no\\
&&+\frac{1}{2}\partial_{x_ix_j}^2F(x,\mu)(\sigma_1\sigma_1^*)^{ij}(s,x,\mu)+\int_{\mR^n}(\partial_\mu F)_i(x,\mu)(u)b_1^i(s,u,\mu)\mu(\dif u)\no\\
&&+\frac{1}{2}\int_{\mR^n}\partial_{u_i}(\partial_\mu F)_j(x,\mu)(u)(\sigma_0\sigma_0^*)^{ij}(s,u,\mu)\mu(\dif u)\no\\
&&+\frac{1}{2}\int_{\mR^n}\partial_{u_i}(\partial_\mu F)_j(x,\mu)(u)(\sigma_1\sigma_1^*)^{ij}(s,u,\mu)\mu(\dif u).
\label{mldefi}
\ee
\et
\begin{proof}
By the extended It\^{o}'s formula in \cite[Proposition 2.9]{dq2}, we know that for $F\in \sS(\mR^n\times\cP_2(\mR^n))$
\be
F(X_t,\sL^{\mP}_{X_t})&=&F(X_0,\sL^{\mP}_{X_0})+\int_0^t(\mL_s F)(X_s,\sL^{\mP}_{X_s})\dif s\no\\
&&+\int_0^t\partial_{x_i} F(X_s,\sL^{\mP}_{X_s})\sigma_0^{ij}(s,X_s,\sL^{\mP}_{X_s})\dif W_s^j\no\\
&&+\int_0^t\partial_{x_i} F(X_s,\sL^{\mP}_{X_s})\sigma_1^{ik}(s,X_s,\sL^{\mP}_{X_s})\dif V_s^k\no\\
&=:&F(X_0,\sL^{\mP}_{X_0})+\int_0^t(\mL_s F)(X_s,\sL^{\mP}_{X_s})\dif s+\Pi_t,
\label{ito1}
\ee
where $(\Pi_t)$ is an $(\mathscr{F}_t)$-adapted local martingale. Thus, by taking the conditional expectation with respect to $\mathscr{F}_t^Y$ on two sides of the above equality, it holds that
\ce
\mE[F(X_t,\sL^{\mP}_{X_t})|\mathscr{F}_t^Y]&=&\mE[F(X_0,\sL^{\mP}_{X_0})|\mathscr{F}_t^Y]+\mE\left[\int_0^t(\mL_s F)(X_s,\sL^{\mP}_{X_s})\dif s|\mathscr{F}_t^Y\right]
+\mE[\Pi_t|\mathscr{F}_t^Y].
\de
We rewrite the above equality to furthermore obtain that
\ce
&&\mE[F(X_t,\sL^{\mP}_{X_t})|\mathscr{F}_t^Y]-\mE[F(X_0,\sL^{\mP}_{X_0})|\mathscr{F}_t^Y]-\int_0^t\mE\left[(\mL_s F)(X_s,\sL^{\mP}_{X_s})|\mathscr{F}_s^Y\right]\dif s\\
&=&\mE\left[\int_0^t(\mL_s F)(X_s,\sL^{\mP}_{X_s})\dif s|\mathscr{F}_t^Y\right]-\int_0^t\mE\left[(\mL_s F)(X_s,\sL^{\mP}_{X_s})|\mathscr{F}_s^Y\right]\dif s
+\mE[\Pi_t|\mathscr{F}_t^Y].
\de
Note that  the right hand side of the above equality is an $(\mathscr{F}^Y_t)$-adapted local martingale (c.f. \cite[Lemma 2.4 and 2.5]{q1}). Hence, by Corollary III 4.27 in \cite{jjas} we have that there exists an $m$-dimensional $(\mathscr{F}^Y_t)$-adapted process $(\Phi_t)$  such that 
\ce
\mE[F(X_t,\sL^{\mP}_{X_t})|\mathscr{F}_t^Y]-\mE[F(X_0,\sL^{\mP}_{X_0})|\mathscr{F}_t^Y]-\int_0^t\mE\left[(\mL_s F)(X_s,\sL^{\mP}_{X_s})|\mathscr{F}_s^Y\right]\dif s=\int_0^t\Phi^*_s\dif\bar{V}_s,
\de
and then
\ce
\mE[F(X_t,\sL^{\mP}_{X_t})|\mathscr{F}_t^Y]=\mE[F(X_0,\sL^{\mP}_{X_0})|\mathscr{F}_t^Y]+\int_0^t\mE\left[(\mL_s F)(X_s,\sL^{\mP}_{X_s})|\mathscr{F}_s^Y\right]\dif s+\int_0^t\Phi^*_s\dif\bar{V}_s.
\de
Since $X_0$ is independent of $(\mathscr{F}^Y_t)$, it holds that
\ce
\mE[F(X_0,\sL^{\mP}_{X_0})|\mathscr{F}_t^Y]=\mE[F(X_0,\sL^{\mP}_{X_0})]=\mE[F(X_0,\sL^{\mP}_{X_0})|\mathscr{F}_0^Y].
\de
From this, it follows that 
\be
\mE[F(X_t,\sL^{\mP}_{X_t})|\mathscr{F}_t^Y]=\mE[F(X_0,\sL^{\mP}_{X_0})|\mathscr{F}_0^Y]+\int_0^t\mE\left[(\mL_s F)(X_s,\sL^{\mP}_{X_s})|\mathscr{F}_s^Y\right]\dif s+\int_0^t\Phi^*_s\dif\bar{V}_s.
\label{pun}
\ee

In the following, we determine the process $(\Phi_t)$. On one side, one can apply the It\^o formula to $\mE[F(X_t,\sL^{\mP}_{X_t})|\mathscr{F}_t^Y]\tilde{V}^j_t$ and obtain that 
\be
&& \mE[F(X_t,\sL^{\mP}_{X_t})|\mathscr{F}_t^Y]\tilde{V}^j_t\no\\
 &=&\int_0^t\mE[F(X_s,\sL^{\mP}_{X_s})|\mathscr{F}_s^Y]\dif \tilde{V}^j_s+\int_0^t\tilde{V}^j_s\dif \mE[F(X_s,\sL^{\mP}_{X_s})|\mathscr{F}_s^Y]+\int_0^t\Phi^j_s\dif s\no\\
 &=&\int_0^t\mE[F(X_s,\sL^{\mP}_{X_s})|\mathscr{F}_s^Y]\mE[h^j(s,X_s,\sL^{\mP}_{X_s},Y_s)|\mathscr{F}_s^Y]\dif s\no\\
 &&+\int_0^t\tilde{V}^j_s\mE\left[(\mL_s F)(X_s,\sL^{\mP}_{X_s})|\mathscr{F}_s^Y\right]\dif s+\int_0^t\Phi^j_s\dif s+I^1_t,
\label{cont1}
\ee
where $(I^1_t)$ is an $(\mathscr{F}^Y_t)$-adapted local martingale.

On the other side, by (\ref{tilw}), (\ref{ito1}) and the It\^o formula for $F(X_t,\sL^{\mP}_{X_t})\tilde{V}^j_t$, we get that for $j=1, 2, \cdots, m$,
\ce
&&F(X_t,\sL^{\mP}_{X_t})\tilde{V}^j_t\\
&=&\int_0^tF(X_s,\sL^{\mP}_{X_s})\dif \tilde{V}^j_s+\int_0^t\tilde{V}^j_s\dif F(X_s,\sL^{\mP}_{X_s})+\int_0^t\partial_{x_i} F(X_s,\sL^{\mP}_{X_s})\sigma^{ij}_1(s,X_s,\sL^{\mP}_{X_s})\dif s\\
&=&\int_0^tF(X_s,\sL^{\mP}_{X_s})h^j(s,X_s,\sL^{\mP}_{X_s},Y_s)\dif s+\int_0^t\tilde{V}^j_s(\mL_s F)(X_s,\sL^{\mP}_{X_s})\dif s\\
&&+\int_0^t\partial_{x_i} F(X_s,\sL^{\mP}_{X_s})\sigma^{ij}_1(s,X_s,\sL^{\mP}_{X_s})\dif s+\int_0^tF(X_s,\sL^{\mP}_{X_s})\dif V^j_s+\int_0^t\tilde{V}^j_s\dif \Pi_s.
\de
Note that $\tilde{V}^j_t$ is measurable with respect to $\mathscr{F}_t^Y$. Thus, by taking the conditional expectation with respect to $\mathscr{F}_t^Y$ on two sides of the above equality, it holds that 
\be
\mE[F(X_t,\sL^{\mP}_{X_t})|\mathscr{F}_t^Y]\tilde{V}^j_t&=&\int_0^t\mE\left[F(X_s,\sL^{\mP}_{X_s})h^j(s,X_s,\sL^{\mP}_{X_s},Y_s)|\mathscr{F}_s^Y\right]\dif s\no\\
&&+\int_0^t\tilde{V}^j_s\mE\left[(\mL_s F)(X_s,\sL^{\mP}_{X_s})|\mathscr{F}_s^Y\right]\dif s\no\\
&&+\int_0^t\mE\left[\partial_{x_i} F(X_s,\sL^{\mP}_{X_s})\sigma^{ij}_1(s,X_s,\sL^{\mP}_{X_s})|\mathscr{F}_s^Y\right]\dif s+I^2_t,
\label{cont2}
\ee
where $(I^2_t)$ denotes an $(\mathscr{F}^Y_t)$-adapted local martingale.

Since the left side of (\ref{cont1}) is the same to that of (\ref{cont2}), bounded variation parts of their right sides should be the same. Therefore,
\be
\Phi^j_s&=&\mE\left[F(X_s,\sL^{\mP}_{X_s})h^j(s,X_s,\sL^{\mP}_{X_s},Y_s)|\mathscr{F}_s^Y\right]-\mE[F(X_s,\sL^{\mP}_{X_s})|\mathscr{F}_s^Y]\mE[h^j(s,X_s,\sL^{\mP}_{X_s},Y_s)|\mathscr{F}_s^Y]\no\\
&&+\mE\left[\partial_{x_i} F(X_s,\sL^{\mP}_{X_s})\sigma^{ij}_1(s,X_s,\sL^{\mP}_{X_s})|\mathscr{F}_s^Y\right], \quad a.s. \mP.
\label{cpro}
\ee
Inserting (\ref{cpro}) in (\ref{pun}) and noting $<\Lambda_t,F>=\mE[F(X_t,\sL^{\mP}_{X_t})|\mathscr{F}_t^Y]$, we get (\ref{kseq0}). Thus, the proof is complete.
\end{proof}

\subsection{The space-distribution dependent Zakai equation}\label{zaka}

Set
\ce
<\tilde{\Lambda}_t,F>:=\mE^{\tilde{\mP}}[F(X_t,\sL^{\mP}_{X_t})\Gamma_t|\mathscr{F}_t^Y], \quad F\in\cB_b(\mR^n\times\cP_2(\mR^n)),
\de
where $\mE^{\tilde{\mP}}$ denotes the expectation under the probability measure $\tilde{\mP}$. Then the space-distribution dependent Zakai equation $\tilde{\Lambda}_{\cdot}$ satisfying is presented as follows.

\bt\label{zakait}(The space-distribution dependent Zakai equation)\\
The space-distribution dependent Zakai equation of the system (\ref{Eq1}) is given by
\be
<\tilde{\Lambda}_t,F>&=&<\tilde{\Lambda}_0,F>+\int_0^t<\tilde{\Lambda}_s,\mL_s F>\dif s
+\int_0^t<\tilde{\Lambda}_s,Fh^j(s,\cdot,\cdot,Y_s)>\dif \tilde{V}^j_s\no\\
&&+\int_0^t<\tilde{\Lambda}_s,\partial_{x_i} F\sigma^{ij}_1(s,\cdot,\cdot)>\dif \tilde{V}^j_s, F\in\sS(\mR^n\times\cP_2(\mR^n)), t\in[0,T].
\label{zakaieq0}
\ee
\et
\begin{proof}
Although the deduction of the space-distribution dependent Zakai equation (\ref{zakaieq0}) is the same to that in \cite[Theorem 2.8]{q1}, we give the proof to the readers' convenience.

By the Kallianpur-Striebel formula, it holds that
\ce
<\Lambda_t,F>=\mE[F(X_t,\sL^{\mP}_{X_t})|\mathscr{F}_t^Y]= \frac{\mE^{\tilde{\mP}}[F(X_t,\sL^{\mP}_{X_t})\Gamma_t|\mathscr{F}_t^Y]}{\mE^{\tilde{\mP}}[\Gamma_t|\mathscr{F}_t^Y]}
= \frac{<\tilde{\Lambda}_t,F>}{<\tilde{\Lambda}_t,1>},
\de
and $<\tilde{\Lambda}_t,F>=<\Lambda_t,F><\tilde{\Lambda}_t,1>$. Thus, to establish the space-distribution dependent Zakai equation (\ref{zakaieq0}), we investigate $<\tilde{\Lambda}_t,1>$. 

First of all, by the It\^o formula, it holds that 
\ce
\Gamma_t=1+\int_0^t \Gamma_sh^i(s,X_s,\sL^{\mP}_{X_s},Y_s)\dif \tilde{V}^i_s.
\de
Taking the conditional expectation with respect to $\sF_t^Y$ under the probability measure $\tilde{\mP}$, one can have that
\ce
\mE^{\tilde{\mP}}[\Gamma_t | \sF_t^Y] = 1 + \int_0^t\mE^{\tilde{\mP}}[\Gamma_sh^i(s,X_s,\sL^{\mP}_{X_s},Y_s) | \sF_s^Y] \dif \tilde{V}^i_s,
\de
namely,
\be
<\tilde{\Lambda}_t,1>= 1 + \int_0^t <\tilde{\Lambda}_s,1><\Lambda_s,h^i(s,\cdot,\cdot,Y_s)>\dif \tilde{V}^i_s.
\label{lamexp}
\ee

Next, combining (\ref{kseq0}) and (\ref{lamexp}) and applying the It\^o formula to $<\Lambda_t,F><\tilde{\Lambda}_t,1>$, we obtain that
 \ce
 &&<\Lambda_t,F><\tilde{\Lambda}_t,1>\\
 &=&<\Lambda_0,F><\tilde{\Lambda}_0,1>+ \int_0^t <\Lambda_s,F>\dif <\tilde{\Lambda}_s,1>+ \int_0^t <\tilde{\Lambda}_s,1>\dif <\Lambda_s,F>\\
 &&+ \int_0^t <\tilde{\Lambda}_s,1><\Lambda_s,h^i(s,\cdot,\cdot,Y_s)>\Phi^i_s\dif s\\
&=&<\Lambda_0,F><\tilde{\Lambda}_0,1>+ \int_0^t <\Lambda_s,F><\tilde{\Lambda}_s,1><\Lambda_s,h^i(s,\cdot,\cdot,Y_s)>\dif \tilde{V}^i_s\\
&&+\int_0^t<\tilde{\Lambda}_s,1><\Lambda_s,\mL_s F>\dif s+\int_0^t<\tilde{\Lambda}_s,1>\Phi^j_s\dif \bar{V}^j_s\\
&&+ \int_0^t <\tilde{\Lambda}_s,1><\Lambda_s,h^i(s,\cdot,\cdot,Y_s)>\Phi^i_s\dif s\\
 &=&<\Lambda_0,F><\tilde{\Lambda}_0,1>+\int_0^t<\tilde{\Lambda}_s,1><\Lambda_s,\mL_s F>\dif s\\
 &&+\int_0^t<\tilde{\Lambda}_s,1><\Lambda_s,\partial_{x_i} F\sigma^{ij}_1(s,\cdot,\cdot)>\dif \tilde{V}^j_s\\
 &&+\int_0^t<\tilde{\Lambda}_s,1><\Lambda_s,Fh^j(s,\cdot,\cdot,Y_s)>\dif \tilde{V}^j_s,
 \de
 where $\Phi$ is defined in (\ref{cpro}). Thus, the above equality together with $<\tilde{\Lambda}_t,F>=<\Lambda_t,F><\tilde{\Lambda}_t,1>$ yields that
\ce
<\tilde{\Lambda}_t,F>&=&<\tilde{\Lambda}_0,F>+ \int_0^t <\tilde{\Lambda}_s,\mL_s F>\dif s+\int_0^t <\tilde{\Lambda}_s,\partial_{x_i} F\sigma^{ij}_1(s,\cdot,\cdot)>\dif \tilde{V}^j_s\\
&&+\int_0^t <\tilde{\Lambda}_s,Fh^j(s,\cdot,\cdot,Y_s)>\dif \tilde{V}^j_s,
\de
 which is just the Zakai equation (\ref{zakaieq0}). The proof is complete.
\end{proof}

\br
If $\sigma_1=0$, Eq.(\ref{kseq0}), (\ref{zakaieq0}) are similar to Eq.(4.36), (4.31) in \cite{sc}, respectively. If $b_1, \sigma_0, \sigma_1, b_2$ are independent of the distribution of $X_t$, Eq.(\ref{kseq0}), (\ref{zakaieq0}) are the same to Eq.(6), (15) without jumps in \cite{q1}, respectively. Therefore, our results are more general.
\er

\section{The pathwise uniqueness for strong solutions to the space-distribution dependent Kushner-Stratonovich equation and the space-distribution dependent Zakai equation}\label{unzaks}

In the section we require that $b_2(t,x,\mu,y), \sigma_2(t,y)$ are independent of $y$. That is, $h(t,x,\mu,y)=h(t,x,\mu)$. First of all, we define strong solutions of the space-distribution dependent Kushner-Stratonovich equation and the space-distribution dependent Zakai equation. After this, a filtered martingale problem is introduced and applied to show the pathwise uniqueness for strong solutions to the space-distribution dependent Kushner-Stratonovich equation and the space-distribution dependent Zakai equation.

\bd\label{ksstrosolu}
A strong solution for the space-distribution dependent Kushner-Stratonovich equation (\ref{kseq0}) is a $(\mathscr{F}_t^Y)_{t\in[0,T]}$-adapted, continuous and $\cP(\mR^n\times\cP_2(\mR^n))$-valued process $(\Pi_t)_{t\in[0,T]}$ such that $(\Pi_t)_{t\in[0,T]}$ solves
the space-distribution dependent Kushner-Stratonovich equation (\ref{kseq0}), that is,
\be
<\Pi_t,F>&=&<\Lambda_0,F>+\int_0^t<\Pi_s,\mL_s F>\dif s+\int_0^t<\Pi_s,\partial_{x_i} F\sigma^{ij}_1(s,\cdot,\cdot)>\dif  \hat{V}^j_s\no\\
&&+\int_0^t\left(<\Pi_s,Fh^j(s,\cdot,\cdot)>-<\Pi_s,F><\Pi_s,h^j(s,\cdot,\cdot)>\right)\dif \hat{V}^j_s, \no\\
&&F\in\sS(\mR^n\times\cP_2(\mR^n)), \qquad t\in[0,T], 
\label{kseqstrsolueq}
\ee
where $\hat{V}_t:=\tilde{V}_t-\int_0^t<\Pi_s,h(s,\cdot,\cdot)>\dif s$.
\ed

\br\label{ksyistrosolu}
By the deduction in Section \ref{nonfilter}, it is obvious that $(\Lambda_t)_{t\in[0,T]}$ is a strong solution of the space-distribution dependent Kushner-Stratonovich equation (\ref{kseq0}).
\er

\bd\label{zakastrosolu}
A strong solution for the space-distribution dependent Zakai equation (\ref{zakaieq0}) is a $(\mathscr{F}_t^Y)_{t\in[0,T]}$-adapted, continuous and $\cM(\mR^n\times\cP_2(\mR^n))$-valued process $(\Sigma_t)_{t\in[0,T]}$ such that $(\Sigma_t)_{t\in[0,T]}$
solves the space-distribution dependent Zakai equation (\ref{zakaieq0}), that is,
\be
<\Sigma_t,F>&=&<\tilde{\Lambda}_0,F>+\int_0^t<\Sigma_s,\mL_s F>\dif s
+\int_0^t<\Sigma_s,Fh^j(s,\cdot,\cdot)>\dif \tilde{V}^j_s\no\\
&&+\int_0^t<\Sigma_s,\partial_{x_i} F\sigma^{ij}_1(s,\cdot,\cdot)>\dif \tilde{V}^j_s, F\in\sS(\mR^n\times\cP_2(\mR^n)), t\in[0,T].
\label{zakastrosolueq}
\ee
\ed

\br\label{zakaiyistrosolu}
By the deduction in Section \ref{nonfilter}, it is obvious that $(\tilde{\Lambda}_t)_{t\in[0,T]}$ is a strong solution of the space-distribution dependent Zakai equation (\ref{zakaieq0}).
\er

Next, we introduce an operator used in the sequel. For $\Phi\in C^{2,2,2}(\mR^n\times\cP_2(\mR^n)\times\mR^m)$, define
\ce
(\cL^{X,\sL^{\mP}_{X},Y}\Phi)(x,\mu,y)&:=&\partial_{x_i}\Phi(x,\mu,y)b^i_1(s,x,\mu)+\frac{1}{2}\partial_{x_ix_j}^2\Phi(x,\mu,y)(\sigma_0\sigma_0^*)^{ij}(s,x,\mu)\\
&&+\frac{1}{2}\partial_{x_ix_j}^2\Phi(x,\mu,y)(\sigma_1\sigma_1^*)^{ij}(s,x,\mu)\\
&&+\int_{\mR^n}(\partial_\mu \Phi)_i(x,\mu,y)(u)b_1^i(s,u,\mu)\mu(\dif u)\\
&&+\frac{1}{2}\int_{\mR^n}\partial_{u_i}(\partial_\mu \Phi)_j(x,\mu,y)(u)(\sigma_0\sigma_0^*)^{ij}(s,u,\mu)\mu(\dif u)\\
&&+\frac{1}{2}\int_{\mR^n}\partial_{u_i}(\partial_\mu \Phi)_j(x,\mu,y)(u)(\sigma_1\sigma_1^*)^{ij}(s,u,\mu)\mu(\dif u)\\
&&+\partial_{y_l}\Phi(x,\mu,y)b^l_2(s,x,\mu)+\frac{1}{2}\partial_{y_ly_q}^2\Phi(x,\mu,y)(\sigma_2\sigma^*_2)^{lq}(s)\\
&&+\partial_{x_iy_l}^2\Phi(x,\mu,y)\sigma^{ij}_1(s,x,\mu)\sigma^{lj}_2(s),
\de
and then we present the concept of filtered martingale problems with respect to $\cL^{X,\sL^{\mP}_{X},Y}$.

\bd\label{fmp}
A process $(\hat{\Pi},\hat{U})$ defined on a probability space $(\hat{\Omega}, \hat{\sF}, (\hat{\sF}_t)_{t\in[0,T]},\hat{\mP})$, with continuous trajectories and values in $\cP(\mR^n\times\cP_2(\mR^n))\times\mR^m$, is a solution of the filtered martingale problem (FMP for short) $(\cL^{X,\sL^{\mP}_{X},Y},X_0,\sL^{\mP}_{X_0},Y_0)$ if 

(i) $\hat{\Pi}$ is $(\sF_t^{\hat{U}})_{t\in[0,T]}$-adapted, 

(ii) for all $\Phi\in\cD(\cL^{X,\sL^{\mP}_{X},Y})$,
\ce
<\hat{\Pi}_t, \Phi\left(\cdot,\cdot,\hat{U}_t\right)>-\int_0^t<\hat{\Pi}_s,\cL^{X,\sL^{\mP}_{X},Y}\Phi(\cdot,\cdot,\hat{U}_s)>\dif s
\de
is a $(\hat{\mP},(\sF_t^{\hat{U}})_{t\in[0,T]})$-martingale,

(iii) for all $\Phi\in\cB_b(\mR^n\times\cP_2(\mR^n)\times\mR^m)$, $\mE^{\hat{\mP}}[<\hat{\Pi}_0,\Phi(\cdot,\cdot,\hat{U}_0)>]=\mE[\Phi(X_0,\sL^{\mP}_{X_0},Y_0)]$.
\ed

\br\label{filtmartprobsolu}
By the deduction in Section \ref{nonfilter}, we see that $(\Lambda, Y)$ is a solution of the FMP $(\cL^{X,\sL^{\mP}_{X},Y},X_0,\sL^{\mP}_{X_0},Y_0)$.
\er

\bd\label{fmpuni}
The uniqueness for the FMP $(\cL^{X,\sL^{\mP}_{X},Y},X_0,\sL^{\mP}_{X_0},Y_0)$ means that if $(\hat{\Pi}^1, \hat{U}^1),\\ (\hat{\Pi}^2, \hat{U}^2)$ defined on these probability spaces $(\hat{\Omega}^1, \hat{\sF}^1, (\hat{\sF}^1_t)_{t\in[0,T]},\hat{\mP}^1)$, $(\hat{\Omega}^2, \hat{\sF}^2, (\hat{\sF}^2_t)_{t\in[0,T]},\hat{\mP}^2)$, respectively, are two solutions of the FMP $(\cL^{X,\sL^{\mP}_{X},Y},X_0,\sL^{\mP}_{X_0},Y_0)$, for $t\in[0,T]$ there exists a Borel measurable 
$\Psi_t: C([0,T],\mR^m)\mapsto\cP(\mR^n\times\cP_2(\mR^n))$ such that 
$$
\hat{\Pi}^1=\Psi_t(U^1), \quad \hat{\mP}^1-a.s., \quad \hat{\Pi}^2=\Psi_t(U^2), \quad \hat{\mP}^2-a.s..
$$
\ed

Here we state and prove the pathwise uniqueness of strong solutions for the space-distribution dependent Kushner-Stratonovich equations by means of filtered martingale problems.

\bt\label{ksunique}
Suppose that the uniqueness holds for the FMP $(\cL^{X,\sL^{\mP}_{X},Y},X_0,\sL^{\mP}_{X_0},Y_0)$. If $\{\Pi_t\}_{t\in[0,T]}$ is a strong solution
of the space-distribution dependent Kushner-Stratonovich equation (\ref{kseq0}). Then $\Pi_t=\Lambda_t, \mP$-a.s. for all $t\in[0,T]$.
\et
\begin{proof}
First of all, for $Y_{\cdot}$ in the system (\ref{Eq1}), i.e.
$$
Y_t=Y_0+\int_0^t b_2(s,X_s,\sL^{\mP}_{X_s})\dif s+\int_0^t\sigma_2(s)\dif V_s,
$$
we apply the It\^o formula to $G(Y_t)$ for any $G\in C_c^\infty(\mR^m)$, and obtain that
\ce
G(Y_t)&=&G(Y_0)+\int_0^t\partial_{y_i} G(Y_s)b^i_2(s,X_s,\sL^{\mP}_{X_s})\dif s
+\frac{1}{2}\int_0^t\partial_{y_iy_j}^2G(Y_s)(\sigma_2\sigma^*_2)^{ij}(s)\dif s\\
&&+\int_0^t\partial_{y_i}G(Y_s)\sigma^{ij}_2(s)\dif V^j_s.
\de
Besides, note that $\{\Pi_t\}_{t\in[0,T]}$ satisfies Eq.(\ref{kseqstrsolueq}), i.e.
\ce
<\Pi_t,F>&=&<\Lambda_0,F>+\int_0^t<\Pi_s,\mL_s F>\dif s+\int_0^t<\Pi_s,\partial_{x_i} F\sigma^{ij}_1(s,\cdot,\cdot)>\dif  \hat{V}^j_s\no\\
&&+\int_0^t\left(<\Pi_s,Fh^j(s,\cdot,\cdot)>-<\Pi_s,F><\Pi_s,h^j(s,\cdot,\cdot)>\right)\dif \hat{V}^j_s, \no\\
&&F\in\sS(\mR^n\times\cP_2(\mR^n)), \qquad t\in[0,T].
\de
Thus, it follows from the It\^o formula that
\be
<\Pi_t,F>G(Y_t)&=&<\Lambda_0,F>G(Y_0)+\int_0^t<\Pi_s,\cL^{X,\sL^{\mP}_{X},Y}(F(\cdot,\cdot)G(Y_s))>\dif s\no\\
&&+\int_0^t<\Pi_s,F>\partial_{y_i}G(Y_s)\sigma^{ij}_2(s)\dif \hat{V}^j_s+\int_0^tG(Y_s)<\Pi_s,\partial_{x_i} F\sigma^{ij}_1(s,\cdot,\cdot)>\dif  \hat{V}^j_s\no\\
&&+\int_0^tG(Y_s)\left[<\Pi_s,Fh^j(s,\cdot,\cdot)>-<\Pi_s,F><\Pi_s,h^j(s,\cdot,\cdot)>\right]\dif \hat{V}^j_s.
\label{pro}
\ee

Now, we observe that
\ce
\hat{V}_t&=&\tilde{V}_t-\int_0^t<\Pi_s,h(s,\cdot,\cdot)>\dif s\\
&=&V_t+\int_0^t h(s,X_s,\sL^{\mP}_{X_s})\dif s-\int_0^t<\Pi_s,h(s,\cdot,\cdot)>\dif s\\
&=&\bar{V}_t-\int_0^t\left(<\Pi_s, h(s,\cdot,\cdot)>-<\Lambda_s, h(s,\cdot,\cdot)>\right)\dif s.
\de
Set
\ce
&&g(s):=<\Pi_s, h(s,\cdot,\cdot)>-<\Lambda_s, h(s,\cdot,\cdot)>,\\
&&\tau_N:=T\wedge\inf\left\{t>0:\int_0^t\left|g(s)\right|^2\dif s>N\right\},
\de
and then $\tau_N$ is a $(\mathscr{F}_t^Y)_{t\in[0,T]}$-stopping time and $\tau_N\rightarrow T$ as $N\rightarrow\infty$
by ($\mathbf{H}^2_{b_2, \sigma_2}$). Define the probability measure
\ce
\frac{\dif\mP_N}{\dif\mP}=\exp\left\{\int_0^{\tau_N}g(s)\dif \bar{V}_s-\frac{1}{2}\int_0^{\tau_N}\left|g(s)\right|^2\dif s\right\}.
\de
Thus, by Lemma \ref{brmoposs} and the Girsanov theorem, it holds that $\hat{V}_t$ is a $(\mathscr{F}_t^Y)_{t\in[0,T]}$-Brownian motion
under $\mP_N$.

Next, for (\ref{pro}) we know that under $\mP_N$,
\ce
<\Pi_{\tau_N\wedge t},F>G(Y_{\tau_N\wedge t})-\int_0^{\tau_N\wedge t}<\Pi_s,\cL^{X,\sL^{\mP}_{X},Y}(F(\cdot,\cdot)G(Y_s))>\dif s
\de
is a $(\mathscr{F}_t^Y)_{t\in[0,T]}$-martingale. Thus, by the appropriate approximation it holds that for $\Phi\in\cD(\cL^{X,\sL^{\mP}_{X},Y})$,
\ce
<\Pi_{\tau_N\wedge t}, \Phi\left(\cdot,\cdot,Y_{\tau_N\wedge t}\right)>-\int_0^{\tau_N\wedge t}<\Pi_s,\cL^{X,\sL^{\mP}_{X},Y}\Phi(\cdot,\cdot,Y_s)>\dif s
\de
is a $(\mP_N, (\mathscr{F}_t^Y)_{t\in[0,T]})$-martingale. Therefore, $(\Pi, Y)$ is a solution of the FMP $(\cL^{X,\sL^{\mP}_{X},Y},X_0,\\\sL^{\mP}_{X_0},Y_0)$ on the probability space $(\Omega,\sF,\{\sF_t\}_{t\in[0,T]},\mP_N)$. Besides, by Remark \ref{filtmartprobsolu}, we see that $(\Lambda, Y)$ is also a solution of the FMP $(\cL^{X,\sL^{\mP}_{X},Y},X_0,\sL^{\mP}_{X_0},Y_0)$ on the probability space $(\Omega,\sF,\{\sF_t\}_{t\in[0,T]},\mP)$. So, by Corollary 3.4 in \cite{ko} and uniqueness 
for the FMP $(\cL^{X,\sL^{\mP}_{X},Y},X_0,\sL^{\mP}_{X_0},Y_0)$, there exists a Borel measurable 
$\Psi_t: C([0,T],\mR^m)\mapsto\cP(\mR^n\times\cP_2(\mR^n))$ such that 
$$
\Pi_t1_{\{t<\tau_N\}}=\Psi_t(Y)1_{\{t<\tau_N\}}, \quad \mP_N-a.s., \quad \Lambda_t=\Psi_t(Y), \quad \mP-a.s.,
$$
and furthermore by the equivalence of $\mP_N$ and $\mP$
$$
\Pi_t1_{\{t<\tau_N\}}=\Lambda_t1_{\{t<\tau_N\}}, \qquad \mP-a.s..
$$
Taking the limits on two sides as $N\rightarrow\infty$, we have
$$
\Pi_t=\Lambda_t, \qquad \mP-a.s..
$$
The proof is complete.
\end{proof}

In the following, we prove the pathwise uniqueness of strong solutions for the space-distribution dependent Zakai equations by means of the above theorem.

\bt\label{unza}
Suppose that the uniqueness holds for the FMP $(\cL^{X,\sL^{\mP}_{X},Y},X_0,\sL^{\mP}_{X_0},Y_0)$. Let $\{\Sigma_t\}_{t\in[0,T]}$ be a strong solution
of the space-distribution dependent Zakai equation (\ref{zakaieq0}). Then $\Sigma_t=\tilde{\Lambda}_t$, $\tilde{\mP}-a.s.$ for all $t\in[0,T]$.
\et
\begin{proof}
First of all, since $\{\Sigma_t\}_{t\in[0,T]}$ is a strong solution of the space-distribution dependent Zakai equation (\ref{zakaieq0}), by Definition \ref{zakastrosolu} it holds that 
\ce
<\Sigma_t,F>&=&<\tilde{\Lambda}_0,F>+\int_0^t<\Sigma_s,\mL_s F>\dif s
+\int_0^t<\Sigma_s,Fh^j(s,\cdot,\cdot)>\dif \tilde{V}^j_s\no\\
&&+\int_0^t<\Sigma_s,\partial_{x_i} F\sigma^{ij}_1(s,\cdot,\cdot)>\dif \tilde{V}^j_s, \quad F\in\sS(\mR^n\times\cP_2(\mR^n)), t\in[0,T],
\de
and 
\ce
<\Sigma_t,1>=1+\int_0^t<\Sigma_s,h^j(s,\cdot,\cdot)>\dif \tilde{V}^j_s.
\de

Next, for $\varepsilon>0$, define the stopping time
\ce
\tau_\varepsilon:=\inf\{t>0: <\Sigma_t,1><\varepsilon\}\wedge T,
\de
and then $<\Sigma_{t\wedge\tau_\varepsilon},1>\geq \e$. And set
\ce
<\Pi_{t\wedge\tau_\varepsilon}, F>:=\frac{<\Sigma_{t\wedge\tau_\varepsilon}, F>}{<\Sigma_{t\wedge\tau_\varepsilon},1>},  \quad F\in\sS(\mR^n\times\cP_2(\mR^n)),
\de
and then $\Pi_{t\wedge\tau_\varepsilon}$ is a $(\mathscr{F}_t^Y)$-adapted $\cP(\mR^n)$-valued process. Applying the It\^o formula to $\frac{<\Sigma_{t\wedge\tau_\varepsilon}, F>}{<\Sigma_{t\wedge\tau_\varepsilon},1>}$, we obtain that
\ce
<\Pi_{t\wedge\tau_\varepsilon},F>&=&<\Lambda_0,F>+\int_0^{t\wedge\tau_\varepsilon}<\Pi_s,\mL_s F>\dif s+\int_0^{t\wedge\tau_\varepsilon}<\Pi_s,\partial_{x_i} F\sigma^{ij}_1(s,\cdot,\cdot)>\dif  \hat{V}^j_s\no\\
&&+\int_0^{t\wedge\tau_\varepsilon}\left(<\Pi_s,Fh^j(s,\cdot,\cdot)>-<\Pi_s,F><\Pi_s,h^j(s,\cdot,\cdot)>\right)\dif \hat{V}^j_s.
\de
From this, it follows that $\{\Pi_{t\wedge\tau_\varepsilon}\}_{t\in[0,T]}$ is a strong solution of the space-distribution dependent Kushner-Stratonovich equation (\ref{kseq0}). By
Theorem \ref{ksunique}, we have that
\be
\Pi_t1_{\{t<\tau_\varepsilon\}}=\Lambda_t1_{\{t<\tau_\varepsilon\}}, \quad \mP-a.s..
\label{ksequiva}
\ee

In the following, we compare $<\Sigma_{t\wedge\tau_\varepsilon},1>$ with $<\tilde{\Lambda}_{t\wedge\tau_\varepsilon},1>$. First of all, by (\ref{lamexp}), it holds that
\ce
<\tilde{\Lambda}_t,1>=\exp\left\{\int_0^t<\Lambda_s,h^j(s, \cdot, \cdot)>\dif \tilde{V}^j_s-\frac{1}{2}\sum_{j=1}^m\int_0^t<\Lambda_s,h^j(s, \cdot, \cdot)>^2\dif s\right\}.
\de
Therefore, $<\tilde{\Lambda}_t,1>>0$. And applying the It\^o formula to $\frac{<\Sigma_{t\wedge\tau_\varepsilon}, 1>}{<\tilde{\Lambda}_{t\wedge\tau_\varepsilon},1>}$, 
we know that 
\ce
\frac{<\Sigma_{t\wedge\tau_\varepsilon}, 1>}{<\tilde{\Lambda}_{t\wedge\tau_\varepsilon},1>}=1
+\int_0^{t\wedge\tau_\varepsilon}\left[\frac{<\Sigma_s, h^j(s, \cdot, \cdot)>}{<\tilde{\Lambda}_s,1>}
-\frac{<\Sigma_s, 1>}{<\tilde{\Lambda}_s,1>}<\Lambda_s,h^j(s, \cdot, \cdot)>\right]\dif \bar{V}^j_s.
\de
Note that for $t<\tau_\varepsilon$,
\ce
&&\frac{<\Sigma_s, h^j(s, \cdot, \cdot)>}{<\tilde{\Lambda}_s,1>}
-\frac{<\Sigma_s, 1>}{<\tilde{\Lambda}_s,1>}<\Lambda_s,h^j(s, \cdot, \cdot)>\\
&=&\frac{<\Sigma_s, 1>}{<\tilde{\Lambda}_s,1>}<\Pi_s,h^j(s, \cdot, \cdot)>-\frac{<\Sigma_s, 1>}{<\tilde{\Lambda}_s,1>}<\Lambda_s,h^j(s, \cdot, \cdot)>\\
&=&0,
\de
where (\ref{ksequiva}) is used in the last equality. So,
\be
\frac{<\Sigma_{t\wedge\tau_\varepsilon}, 1>}{<\tilde{\Lambda}_{t\wedge\tau_\varepsilon},1>}=1, \quad \mP-a.s..
\label{1equ}
\ee

Based on (\ref{ksequiva}) and (\ref{1equ}), it holds that
\ce
\Sigma_t1_{\{t<\tau_\varepsilon\}}=<\Sigma_t, 1>\Pi_t1_{\{t<\tau_\varepsilon\}}=<\tilde{\Lambda}_t,1>\Lambda_t1_{\{t<\tau_\varepsilon\}}=\tilde{\Lambda}_t1_{\{t<\tau_\varepsilon\}}, \quad \mP-a.s.,
\de
and furthermore by the equivalence between $\tilde{\mP}$ and $\mP$,
$$
\Sigma_t1_{\{t<\tau_\varepsilon\}}=\tilde{\Lambda}_t1_{\{t<\tau_\varepsilon\}}, \quad \tilde{\mP}-a.s..
$$
From this, it follows that $\tau_\varepsilon\geq\inf\{t>0: <\tilde{\Lambda}_s,1><\varepsilon\}\wedge T$. Note that
$\inf\{t>0: <\tilde{\Lambda}_s,1><\varepsilon\}\wedge T=T$ when $\varepsilon$ is small enough. So,
$\tau_\varepsilon=T$ and $\Sigma_t=\tilde{\Lambda}_t, \tilde{\mP}-a.s.$. The proof is complete.
\end{proof}

\br\label{uniqcond}
About the conditions under which the uniqueness holds for the FMP $(\cL^{X,\sL^{\mP}_{X},Y},X_0,\sL^{\mP}_{X_0},Y_0)$, please refer to Theorem 3.2 in \cite{ko}.
\er

\section{Superposition between the space-distribution dependent Zakai equation and a space-distribution dependent Fokker-Planck equation}\label{super}

In the section we also require that $b_2(t,x,\mu,y), \sigma_2(t,y)$ are independent of $y$. First of all, we define weak solutions of the space-distribution dependent Zakai equation (\ref{zakaieq0}). Then a type of space-distribution dependent Fokker-Planck equations and their weak solutions are introduced. Finally, we prove a superposition principle between weak solutions of the space-distribution dependent Zakai equation (\ref{zakaieq0}) and weak solutions of a space-distribution dependent Fokker-Planck equation.

\subsection{The space-distribution dependent Zakai equation}

In the subsection, we view the space-distribution dependent Zakai equation (\ref{zakaieq0}) as a SDE and define its weak solutions and the uniqueness in law.

\bd\label{soluzakai}
$\{(\hat{\Omega}, \hat{\mathscr{F}}, \{\hat{\mathscr{F}}_t\}_{t\in[0,T]},\hat{\mP}), (\hat{\Sigma}_t,
\hat{V}_t)\}$ is called a weak solution of the space-distribution dependent Zakai equation (\ref{zakaieq0}), if the following holds:

(i) $(\hat{\Omega}, \hat{\mathscr{F}}, \{\hat{\mathscr{F}}_t\}_{t\in[0,T]},\hat{\mP})$ is a complete filtered
probability space;

(ii) $\hat{\Sigma}_t$ is a $\cM(\mR^n\times\cP_2(\mR^n))$-valued $(\hat{\mathscr{F}}_t)$-adapted continuous process and $\hat{\Sigma}_0\in\cP(\mR^n\times\cP_2(\mR^n))$;

(iii) $\hat{V}_t$ is a $m$-dimensional $(\hat{\mathscr{F}}_t)$-adapted Brownian motion;

(iv) For any $t\in[0,T]$,
\ce
&&\hat{\mP}\bigg(\int_0^t\int_{\mR^n\times\cP_2(\mR^n)}\Big(|b_1(r,x,\mu)|+|h(r,x,\mu)|^2+\|\sigma_1(r,x,\mu)\|^2\\
&&+\|\sigma_0\sigma_0^*(r,x,\mu)\|+\int_{\mR^n}|b_1(r,u,\mu)|\mu(\dif u)+\int_{\mR^n}\|\sigma_0\sigma_0^*(r,u,\mu)\|\mu(\dif u)\\
&&+\int_{\mR^n}\|\sigma_1\sigma_1^*(r,u,\mu)\|\mu(\dif u)\Big)\hat{\Sigma}_r(\dif x\dif \mu)\dif r<\infty\bigg)=1;
\de

(v) $(\hat{\Sigma}_t, \hat{V}_t)$ satisfies the following equation
\be
<\hat{\Sigma}_t, F>&=&<\hat{\Sigma}_0, F>+\int_0^t<\hat{\Sigma}_s, (\mL_s F)(\cdot,\cdot)>\dif s
+\int_0^t<\hat{\Sigma}_s, F h^j(s,\cdot,\cdot)>\dif \hat{V}^j_s\no\\
&&+\int_0^t<\hat{\Sigma}_s, \partial_{x_i} F\sigma^{ij}_1(s,\cdot,\cdot)>\dif \hat{V}^j_s, t\in[0,T], F\in\sS(\mR^n\times\cP_2(\mR^n)).
\label{zakaieqweaksolu}
\ee
\ed

\br\label{zakaiweak}
By the deduction in Section \ref{nonfilter}, it is obvious that $\{(\Omega, \mathscr{F}, \{\mathscr{F}_t\}_{t\in[0,T]}, \tilde{\mP}),
(\tilde{\Lambda}_t, \tilde{V}_t)\}$ is a weak solution of the space-distribution dependent Zakai equation (\ref{zakaieq0}).
\er

\bd\label{launzakai}
The uniqueness in law of weak solutions for the space-distribution dependent Zakai equation (\ref{zakaieq0}) means that if there exist two weak solutions $\{(\hat{\Omega}^1, \hat{\mathscr{F}}^1, \{\hat{\mathscr{F}}^1_t\}_{t\in[0,T]}, \hat{\mP}^1), \\(\hat{\Sigma}^1_t,\hat{V}^1_t)\}$ and $\{(\hat{\Omega}^2, \hat{\mathscr{F}}^2, \{\hat{\mathscr{F}}^2_t\}_{t\in[0,T]}, \hat{\mP}^2),(\hat{\Sigma}^2_t,\hat{V}^2_t)\}$ with $\hat{\mP}^1\circ(\hat{\Sigma}^1_0)^{-1}=\hat{\mP}^2\circ(\hat{\Sigma}^2_0)^{-1}$, then $\hat{\mP}^1\circ(\hat{\Sigma}^1_{\cdot})^{-1}=\hat{\mP}^2\circ(\hat{\Sigma}^2_{\cdot})^{-1}$.
\ed

\subsection{The space-distribution dependent Fokker-Planck equation}\label{fpes1}

In the subsection, we introduce the space-distribution dependent Fokker-Planck equation associated with the space-distribution dependent Zakai equation (\ref{zakaieq0}) and define its weak solutions.

First of all, set 
\ce
&&\sG:=\Big\{\Sigma\in\cM(\mR^n\times\cP_2(\mR^n))\mapsto G(\Sigma)=g\left(<\Sigma,F_1>,\cdots,<\Sigma,F_k>\right): k\in\mN,\\
&&\qquad\qquad\qquad\qquad\qquad g\in C_b^2(\mR^{k}), F_1, \cdots, F_k\in C_b^{2,2}(\mR^n\times\cP_2(\mR^n))\Big\}.
\de
Then for any $G(\Sigma)=g\left(<\Sigma,F_1>,\cdots,<\Sigma,F_k>\right)=:g(<\Sigma,\boldsymbol{F}>)\in\sG$, we define the operator ${\bf L}_t$ on $\sG$:
\ce
{\bf L}_tG(\Sigma)&=&\frac{1}{2}\partial_{y_uy_v} g(<\Sigma,\boldsymbol{F}>)<\Sigma,F_u h^l(t,\cdot,\cdot)+\partial_{x_i}F_u\sigma^{il}_1(t,\cdot,\cdot)>\\
&&\times<\Sigma,F_v h^l(t,\cdot,\cdot)+\partial_{x_i}F_v\sigma^{il}_1(t,\cdot,\cdot)>\\
&&+\partial_{y_u} g(<\Sigma,\boldsymbol{F}>)<\Sigma,(\mL_tF_u)(\cdot,\cdot)>\\
&\overset{(\ref{mldefi})}{=}&\frac{1}{2}\partial_{y_uy_v} g(<\Sigma,\boldsymbol{F}>)<\Sigma,F_u h^l(t,\cdot,\cdot)+\partial_{x_i}F_u\sigma^{il}_1(t,\cdot,\cdot)>\\
&&\times<\Sigma,F_v h^l(t,\cdot,\cdot)+\partial_{x_i}F_v\sigma^{il}_1(t,\cdot,\cdot)>\\
&&+\partial_{y_u} g(<\Sigma,\boldsymbol{F}>)<\Sigma,\partial_{x_i}F_u b^i_1(t,\cdot,\cdot)>\\
&&+\frac{1}{2}\partial_{y_u} g(<\Sigma,\boldsymbol{F}>)<\Sigma,\partial_{x_ix_j}F_u\left(\sigma_0\sigma_0^*\right)^{ij}(t,\cdot,\cdot)>\\
&&+\frac{1}{2}\partial_{y_u} g(<\Sigma,\boldsymbol{F}>)<\Sigma,\partial_{x_ix_j}F_u\left(\sigma_1\sigma_1^*\right)^{ij}(t,\cdot,\cdot)>\\
&&+\partial_{y_u} g(<\Sigma,\boldsymbol{F}>)<\Sigma,\int_{\mR^n}(\partial_\mu F_u)_i(\cdot,\cdot)(w)b_1^i(s,w,\mu)\mu(\dif w)>\\
&&+\frac{1}{2}\partial_{y_u} g(<\Sigma,\boldsymbol{F}>)<\Sigma,\int_{\mR^n}\partial_{u_i}(\partial_\mu F_u)_j(\cdot,\cdot)(w)(\sigma_0\sigma_0^*)^{ij}(s,w,\mu)\mu(\dif w)>\\
&&+\frac{1}{2}\partial_{y_u} g(<\Sigma,\boldsymbol{F}>)<\Sigma,\int_{\mR^n}\partial_{u_i}(\partial_\mu F_u)_j(\cdot,\cdot)(w)(\sigma_1\sigma_1^*)^{ij}(s,w,\mu)\mu(\dif w)>.
\de
Consider the following Fokker-Planck equation (FPE for short):
\be
\partial_t\Xi_t={\bf L}^*_t\Xi_t,
\label{fpecoup}
\ee
where $(\Xi_t)_{t\in[0,T]}$ is a family of probability measures on $\sB(\cM(\mR^n\times\cP_2(\mR^n)))$. Then, we define weak solutions of Eq.(\ref{fpecoup}).

\bd\label{weaksolufpe}
A measurable family $(\Xi_t)_{t\in[0,T]}$ of probability measures on $\sB(\cM(\mR^n\times\cP_2(\mR^n)))$ is called a weak solution of Eq.(\ref{fpecoup}) starting from $\Xi_0$ at time $0$ if 
\be
&&\int_0^T\int_{\cM(\mR^n\times\cP_2(\mR^n))}\int_{\mR^n\times\cP_2(\mR^n)}\Big(|b_1(r,x,\mu)|+|h(r,x,\mu)|^2+\|\sigma_1(r,x,\mu)\|^2\no\\
&&+\|\sigma_0\sigma_0^*(r,x,\mu)\|+\int_{\mR^n}|b_1(r,u,\mu)|\mu(\dif u)+\int_{\mR^n}\|\sigma_0\sigma_0^*(r,u,\mu)\|\mu(\dif u)\no\\
&&+\int_{\mR^n}\|\sigma_1\sigma_1^*(r,u,\mu)\|\mu(\dif u)\Big)\Sigma(\dif x\dif \mu)\Xi_r(\dif \Sigma)\dif r<\infty, \label{fpem01}
\ee
and for any $G\in\sG$ and $0\leq t\leq T$,
\be
&&\int_{\cM(\mR^n\times\cP_2(\mR^n))}G(\Sigma)\Xi_t(\dif\Sigma)\no\\
&=&\int_{\cM(\mR^n\times\cP_2(\mR^n))}G(\Sigma)\Xi_0(\dif \Sigma)+\int_0^t\int_{\cM(\mR^n\times\cP_2(\mR^n))}{\bf L}_rG(\Sigma)\Xi_r(\dif \Sigma)\dif r.
\label{fpem2}
\ee
The uniqueness of the weak solutions to Eq.(\ref{fpecoup}) means that, if for any $s\in[0,T]$ and any $\Xi_s\in\cP(\cM(\mR^n))$, $(\Xi_t)_{t\in[s,T]}$ and $(\tilde{\Xi}_t)_{t\in[s,T]}$ are two weak solutions to Eq.(\ref{fpecoup}) starting from $\Xi_s$ 
at time $s$, then $\Xi_t=\tilde{\Xi}_t$ for any $t\in[s,T]$.
\ed

It is easy to see that under the condition (\ref{fpem01}), two integrals in the right side of Eq.(\ref{fpem2}) are well defined.

\subsection{A superposition principle between Eq.(\ref{zakaieq0}) and Eq.(\ref{fpecoup})}\label{susupeprin1}

In the subsection, we prove the following superposition principle between Eq.(\ref{zakaieq0}) and Eq.(\ref{fpecoup}). 

\bt(The superposition principle on $\cM(\mR^n\times\cP_2(\mR^n))$)\label{supeprin}

(i) If $\{(\hat{\Omega}, \hat{\mathscr{F}}, \{\hat{\mathscr{F}}_t\}_{t\in[0,T]},\hat{\mP}), (\hat{\Sigma}_t,\hat{V}_t)\}$ is a weak solution  for Eq.(\ref{zakaieq0}), then $(\sL^{\hat{\mP}}_{\hat{\Sigma}_t})$ is a weak solution for Eq.(\ref{fpecoup}).

(ii) If weak solutions for Eq.(\ref{fpecoup}) are unique, then weak solutions  for Eq.(\ref{zakaieq0}) have the uniqueness in law.
\et
\begin{proof}
Since $\{(\hat{\Omega}, \hat{\mathscr{F}}, \{\hat{\mathscr{F}}_t\}_{t\in[0,T]},\hat{\mP}), (\hat{\Sigma}_t,\hat{V}_t)\}$ is a weak solution  for Eq.(\ref{zakaieq0}), by Definition \ref{soluzakai} (iv), it holds that 
\ce
&&\int_0^T\int_{\mR^n\times\cP_2(\mR^n)}\Big(|b_1(r,x,\mu)|+|h(r,x,\mu)|^2+\|\sigma_1(r,x,\mu)\|^2\\
&&+\|\sigma_0\sigma_0^*(r,x,\mu)\|+\int_{\mR^n}|b_1(r,u,\mu)|\mu(\dif u)+\int_{\mR^n}\|\sigma_0\sigma_0^*(r,u,\mu)\|\mu(\dif u)\\
&&+\int_{\mR^n}\|\sigma_1\sigma_1^*(r,u,\mu)\|\mu(\dif u)\Big)\hat{\Sigma}_r(\dif x\dif \mu)\dif r<\infty, a.s.,
\de
and 
\ce
&&\hat{\mE}\int_0^T\int_{\mR^n\times\cP_2(\mR^n)}\Big(|b_1(r,x,\mu)|+|h(r,x,\mu)|^2+\|\sigma_1(r,x,\mu)\|^2\\
&&+\|\sigma_0\sigma_0^*(r,x,\mu)\|+\int_{\mR^n}|b_1(r,u,\mu)|\mu(\dif u)+\int_{\mR^n}\|\sigma_0\sigma_0^*(r,u,\mu)\|\mu(\dif u)\\
&&+\int_{\mR^n}\|\sigma_1\sigma_1^*(r,u,\mu)\|\mu(\dif u)\Big)\hat{\Sigma}_r(\dif x\dif \mu)\dif r\\
&=&\int_0^T\int_{\cM(\mR^n\times\cP_2(\mR^n))}\int_{\mR^n\times\cP_2(\mR^n)}\Big(|b_1(r,x,\mu)|+|h(r,x,\mu)|^2+\|\sigma_1(r,x,\mu)\|^2\\
&&+\|\sigma_0\sigma_0^*(r,x,\mu)\|+\int_{\mR^n}|b_1(r,u,\mu)|\mu(\dif u)+\int_{\mR^n}\|\sigma_0\sigma_0^*(r,u,\mu)\|\mu(\dif u)\\
&&+\int_{\mR^n}\|\sigma_1\sigma_1^*(r,u,\mu)\|\mu(\dif u)\Big)\Sigma(\dif x\dif \mu)\sL^{\hat{\mP}}_{\hat{\Sigma}_r}(\dif \Sigma)\dif r<\infty,
\de
where $\hat{\mE}$ denotes the expectation under the probability $\hat{\mP}$. So, $(\sL^{\hat{\mP}}_{\hat{\Sigma}_t})$ satisfies (\ref{fpem01}).

Next, from Definition \ref{soluzakai} (v), it follows that  for $F\in C_b^{2,2}(\mR^n\times\cP_2(\mR^n))$,
 \ce
<\hat{\Sigma}_t, F>&=&<\hat{\Sigma}_0, F>+\int_0^t<\hat{\Sigma}_s, (\mL_s F)(\cdot,\cdot)>\dif s
+\int_0^t<\hat{\Sigma}_s, F h^j(s,\cdot,\cdot)>\dif \hat{V}^j_s\no\\
&&+\int_0^t<\hat{\Sigma}_s, \partial_{x_i} F\sigma^{ij}_1(s,\cdot,\cdot)>\dif \hat{V}^j_s, \quad t\in[0,T].
\de
For $G(\Sigma)\in\sG$, applying the It\^o formula to $G(\hat{\Sigma}_t)$, we know that
\ce
G(\hat{\Sigma}_t)&=&G(\hat{\Sigma}_0)+\int_0^t\partial_{y_u} g(<\hat{\Sigma}_s,\boldsymbol{F}>)<\hat{\Sigma}_s,(\mL_sF_u)(\cdot,\cdot)>\dif s\\
&&+\int_0^t\partial_{y_u} g(<\hat{\Sigma}_s,\boldsymbol{F}>)<\hat{\Sigma}_s,F_u h^l(s,\cdot,\cdot)+\partial_{x_i}F_u\sigma^{il}_1(s,\cdot,\cdot)>\dif \hat{V}^j_s\\
&&+\frac{1}{2}\int_0^t\partial_{y_uy_v} g(<\hat{\Sigma}_s,\boldsymbol{F}>)<\hat{\Sigma}_s,F_u h^l(s,\cdot,\cdot)+\partial_{x_i}F_u\sigma^{il}_1(s,\cdot,\cdot)>\\
&&\qquad \times<\hat{\Sigma}_s,F_v h^l(s,\cdot,\cdot)+\partial_{x_i}F_v\sigma^{il}_1(s,\cdot,\cdot)>\dif s.
\de
By taking the expectation on two sides, it holds that
\ce
\hat{\mE}G(\hat{\Sigma}_t)&=&\hat{\mE}G(\hat{\Sigma}_0)+\int_0^t\hat{\mE}\partial_{y_u} g(<\hat{\Sigma}_s,\boldsymbol{F}>)<\hat{\Sigma}_s,(\mL_sF_u)(\cdot,\cdot)>\dif s\\
&&+\frac{1}{2}\int_0^t\hat{\mE}\partial_{y_uy_v} g(<\hat{\Sigma}_s,\boldsymbol{F}>)<\hat{\Sigma}_s,F_u h^l(s,\cdot,\cdot)+\partial_{x_i}F_u\sigma^{il}_1(s,\cdot,\cdot)>\\
&&\qquad \times<\hat{\Sigma}_s,F_v h^l(s,\cdot,\cdot)+\partial_{x_i}F_v\sigma^{il}_1(s,\cdot,\cdot)>\dif s,
\de
and furthermore
\ce
&&\int_{\cM(\mR^n\times\cP_2(\mR^n))}G(\Sigma)\sL^{\hat{\mP}}_{\hat{\Sigma}_t}(\dif\Sigma)\no\\
&=&\int_{\cM(\mR^n\times\cP_2(\mR^n))}G(\Sigma)\sL^{\hat{\mP}}_{\hat{\Sigma}_0}(\dif \Sigma)+\int_0^t\int_{\cM(\mR^n\times\cP_2(\mR^n))}{\bf L}_sG(\Sigma)\sL^{\hat{\mP}}_{\hat{\Sigma}_s}(\dif \Sigma)\dif s.
\de
Thus, $(\sL^{\hat{\mP}}_{\hat{\Sigma}_t})$ satisfies (\ref{fpem2}). The proof of (i) is complete.

For (ii), we assume that $\{(\hat{\Omega}^1, \hat{\mathscr{F}}^1, \{\hat{\mathscr{F}}^1_t\}_{t\in[0,T]}, \hat{\mP}^1), (\hat{\Sigma}^1_t,\hat{V}^1_t)\}$ and $\{(\hat{\Omega}^2, \hat{\mathscr{F}}^2, \{\hat{\mathscr{F}}^2_t\}_{t\in[0,T]}, \hat{\mP}^2),\\(\hat{\Sigma}^2_t,\hat{V}^2_t)\}$ are two weak solutions for Eq.(\ref{zakaieq0}) with $\hat{\mP}^1\circ(\hat{\Sigma}^1_0)^{-1}=\hat{\mP}^2\circ(\hat{\Sigma}^2_0)^{-1}$ and $\hat{\mP}^1\circ(\hat{\Sigma}^1_{\cdot})^{-1}\neq\hat{\mP}^2\circ(\hat{\Sigma}^2_{\cdot})^{-1}$. Then it holds that for any $t\in(0,T]$, $\hat{\mP}^1\circ(\hat{\Sigma}^1_t)^{-1}\neq\hat{\mP}^2\circ(\hat{\Sigma}^2_t)^{-1}$. 

Besides, by (i), we know that $\hat{\mP}^1\circ(\hat{\Sigma}^1_t)^{-1}, \hat{\mP}^2\circ(\hat{\Sigma}^2_t)^{-1}$ are two weak solutions for Eq.(\ref{fpecoup}). And the uniqueness of weak solutions for Eq.(\ref{fpecoup}) implies that $\hat{\mP}^1\circ(\hat{\Sigma}^1_t)^{-1}=\hat{\mP}^2\circ(\hat{\Sigma}^2_t)^{-1}$. This contradicts the above conclusion. Therefore, (ii) holds. The proof is complete.
\end{proof}

\br
Since $C_b^{2,2}(\mR^n\times\cP_2(\mR^n))$ is not separable in the usual norms, we can't obtain the converse conclusion.
\er

Combining Remark \ref{zakaiweak} and Theorem \ref{supeprin} (i), we have the following conclusion.

\bc
Under the assumptions ($\mathbf{H}^1_{b_1, \sigma_0, \sigma_1}$) ($\mathbf{H}^2_{b_1, \sigma_0,\sigma_1}$) ($\mathbf{H}^2_{b_2, \sigma_2}$), Eq.(\ref{fpecoup}) has a weak solution.
\ec

\section{Conclusion}\label{con}

In the paper, we consider the space-distribution dependent Zakai equations from nonlinear filtering problems of McKean-Vlasov SDEs with correlated noises. Firstly, we establish the space-distribution dependent Kushner-Stratonovich equations and the space-distribution dependent Zakai equations. Secondly, the pathwise uniqueness of strong solutions for the space-distribution dependent Kushner-Stratonovich equations and the space-distribution dependent Zakai equations is shown. Finally, we prove a superposition principle between the space-distribution dependent Zakai equations and space-distribution dependent Fokker-Planck equations. As a by-product, we give some conditions under which space-distribution dependent Fokker-Planck equations have weak solutions.

Our methods also can be used to solve nonlinear filtering problems of McKean-Vlasov SDEs with correlated sensor noises. Concretely speaking, consider the slow-fast system in the system (\ref{Eq0}), i.e. 
\be\left\{\begin{array}{l}
\dif \check{X}_t=\check{b}_1(t,\check{X}_t, \sL^{\mP}_{\check{X}_t})\dif t+\check{\sigma}_1(t,\check{X}_t, \sL^{\mP}_{\check{X}_t})\dif V_t,\\
\dif \check{Y}_t=\check{b}_2(t,\check{X}_t, \sL^{\mP}_{\check{X}_t},\check{Y}_t)\dif t+\check{\sigma}_2\dif W_t+\check{\sigma}_3\dif V_t, \quad 0\leq t\leq T.
 \end{array}
\right.
\label{0sfs} 
\ee
We assume the following:
\begin{enumerate}[(i)]
\item $\check{b}_1, \check{\sigma}_1$ satisfy ($\mathbf{H}^1_{b_1, \sigma_0, \sigma_1}$)-($\mathbf{H}^2_{b_1,\sigma_0,\sigma_1}$), where $\check{b}_1, \check{\sigma}_1$ replace $b_1, \sigma_1$, respectively;
\item $\check{b}_2(t,x,\mu,y)$ is bounded for all $t\in[0,T], x\in\mR^n, \mu\in \cP_2(\mR^n), y\in\mR^m$;
\item $\check{\sigma}_2\check{\sigma}^*_2+\check{\sigma}_3\check{\sigma}^*_3=I_m,$ where $\check{\sigma}^*_2$ stands for the transpose of the matrix $\check{\sigma}_2$ and $I_m$ is the $m$-order unit matrix.
\end{enumerate}

Under the above assumptions, it holds that the system (\ref{0sfs}) has a unique strong solution denoted as $(\check{X}_t,\check{Y}_t)$. Set
\ce
&&U_t:=\check{\sigma}_2 W_t+\check{\sigma}_3 V_t, \\
&&\check{\Gamma}^{-1}_t:=\exp\bigg\{-\int_0^t\check{b}_2^i(s,\check{X}_s,\sL^{\mP}_{\check{X}_s},\check{Y}_s)\dif U^i_s-\frac{1}{2}\int_0^t
\left|\check{b}_2(s,\check{X}_s,\sL^{\mP}_{\check{X}_s},\check{Y}_s)\right|^2\dif s\bigg\},
\de
and then $\check{\Gamma}^{-1}_{\cdot}$ is an exponential martingale. Moreover, define the probability measure 
$$
\frac{\dif \tilde{\check{\mP}}}{\dif \mP}:=\check{\Gamma}^{-1}_T,
$$
and set
\ce
<\tilde{\check{\Lambda}}_t,F>:=\mE^{\tilde{\check{\mP}}}[F(\check{X}_t,\sL^{\mP}_{\check{X}_t})\check{\Gamma}_t|\mathscr{F}_t^{\check{Y}}],
\de
where $\mE^{\tilde{\check{\mP}}}$ stands for the expectation under the probability measure $\tilde{\check{\mP}}$. By deduction the same as to that in Theorem \ref{zakait}, we obtain the following space-distribution dependent Zakai equation.

\bc (The space-distribution dependent Zakai equation)\label{zasen}\\
 The space-distribution dependent Zakai equation of the system (\ref{Eq0}) is given by
\be
<\tilde{\check{\Lambda}}_t,F>&=&<\tilde{\check{\Lambda}}_0,F>+\int_0^t<\tilde{\check{\Lambda}}_s,\check{\mL}_s F>\dif s+\int_0^t<\tilde{\check{\Lambda}}_s,F(\cdot,\cdot)\check{b}_2^l(s,\cdot,\cdot,\check{Y}_s)>\dif  \tilde{U}^j_s\no\\
&&+\int_0^t<\tilde{\check{\Lambda}}_s,(\partial_{x_i} F)(\cdot,\cdot)\check{\sigma}^{ik}_1(s,\cdot,\cdot)\check{\sigma}^{jk}_2>\dif  \tilde{U}^j_s, \quad t\in[0,T],
\label{zakaieq01}
\ee
where 
\ce
(\check{\mL}_s F)(x,\mu)&:=&\partial_{x_i}F(x,\mu)\check{b}^i_1(s,x,\mu)+\frac{1}{2}\partial_{x_ix_j}^2F(x,\mu)(\check{\sigma}_1\check{\sigma}_1^*)^{ij}(s,x,\mu)\\
&&+\int_{\mR^n}(\partial_\mu F)_i(x,\mu)(y)\check{b}_1^i(s,y,\mu)\mu(\dif y)\\
&&+\frac{1}{2}\int_{\mR^n}\partial_{y_i}(\partial_\mu F)_j(x,\mu)(y)(\check{\sigma}_1\check{\sigma}_1^*)^{ij}(s,y,\mu)\mu(\dif y),
\de
and
$\tilde{U}_t:=U_t+\int_0^t \check{b}_2(s,\check{X}_s,\sL^{\mP}_{\check{X}_s},\check{Y}_s)\dif s$.
\ec

Of course, we can study the pathwise uniqueness and superposition principles for Eq.(\ref{zakaieq01}) by means the same as to that in Theorem \ref{unza} and Theorem \ref{supeprin}.

 \bigskip

\textbf{Acknowledgements:}

The second named author would like to thank Professor Renming Song for providing her an excellent environment to work in the University of Illinois at Urbana-Champaign.

\end{document}